\documentclass[12pt]{amsart}

\begin{document}

\title{Subvarieties of abelian varieties}

\author{E. Izadi}
\address{Department of Mathematics, Boyd
Graduate Studies Research Center, University of Georgia, Athens, GA
30602-7403, USA}
\email{izadi@math.uga.edu}
\thanks{The author was partially supported by grant DMS-0071795 from
the National Science Foundation}
\date{March 15 2001}

%\begin{abstract
%We discuss various constructions which allow one to embed a
%principally polarized abelian variety in the jacobian of a curve. Each
%of these gives representatives of multiples of the minimal cohomology
%class for curves which in turn produce subvarieties of higher
%dimension representing multiples of the minimal class. We then discuss
%the problem of producing curves representing multiples of the minimal
%class via deformation-theoretic methods.
%\end{abstract}

\maketitle

Let $(A ,\Theta)$ be a principally polarized abelian variety ({\em
ppav}) of dimension $g$ over the field ${\mathbb C}$ of complex
numbers. This means that $\Theta$ is an ample divisor on $A$,
well-determined up to translation, with $h^0 (A ,\Theta ) :=
\hbox{dim} H^0 (A ,\Theta) =1$. Let $[\Theta]\in H^2 (A,{\mathbb Z})$
be the cohomology class of the theta divisor $\Theta$. Then the
cohomology class $[\Theta]^{ g-e }$ is divisible by $(g-e)!$. The
class $\frac{[\Theta]^{ g-e}}{ (g-e)!}$ is not divisible and it is
called the dimension $e$ minimal cohomology class in $(A
,\Theta)$. This class is positive in the sense that some multiple of
it can be represented by an algebraic subvariety (for instance
$[\Theta]^{ g-e }$ is the class of a complete intersection of $g-e$
general translates of $\Theta$) and, furthermore, any subvariety whose
class is a multiple of $\frac{[\Theta]^{ g-e }}{ (g-e )!}$ is
nondegenerate, i.e., generates $A$ as a group. We are interested in
the representability of multiples of the minimal classes by algebraic
subvarieties of $A$. We begin by discussing two special cases.

\section{Jacobians}

Let $C$ be a smooth, complete, irreducible curve of genus $g$ over the
complex numbers. The jacobian $JC = Pic^0 C$ of $C$ is the connected
component of its Picard group parametrizing degree $0$ invertible
sheaves. For any nonnegative integer $e$, the choice of an invertible
sheaf ${\mathcal L}$ of degree $e$ on $C$ gives a morphism
\[
\begin{array}{rcl}
\phi_{{\mathcal L}} : C^{ (e) } &\longrightarrow & JC\\
D_e &\longmapsto &{\mathcal O}_C (D_e)\otimes{\mathcal L}^{ -1}
\end{array}
\]
where $C^{ (e) }$ is the $e$-th symmetric power of $C$. For $e\geq g$,
such a morphism is surjective. When $e=g-1$, the image of $C^{ (g-1)
}$ by $\phi_{\mathcal L}$ is a theta divisor on $JC$ which we will
denote by $\Theta_C$ (always well-determined up to translation). For
any $e$ between $1$ and $g$ the image of $\phi_{{\mathcal L}}$ has
class $\frac{[\Theta]^{ g-e}}{ (g-e)!}$. If $e=1$, the map
$\phi_{{\mathcal L}}$ is an embedding and its image is called an Abel
curve. By a theorem of Matsusaka \cite{matsusaka59}, the dimension $1$
minimal class is represented by an algebraic curve $C$ in $(A
,\Theta)$ if and only if $(A ,\Theta)$ is the polarized jacobian $( JC
,\Theta_C )$ of $C$. The higher-dimensional analogue of this theorem
has the following counterexample: by a result of Clemens and
Griffiths, the intermediate jacobian of a smooth cubic threefold in
${\mathbb P}^4$ is not the jacobian of a curve but it contains a
surface (the image of the Fano variety of lines in the cubic
threefold) whose cohomology class is the dimension $2$ minimal class
\cite{clemensgriffiths72}. Debarre has shown that for any $e$ strictly
between $1$ and $g-1$ jacobians form an irreducible component of the
family of ppav in which the dimension $e$ minimal class is represented
by an algebraic subvariety \cite{debarre95}.

As the above suggests, not every ppav is a jacobian. In fact, the
moduli space of ppav of dimension $g$ has dimension $\frac{ g(g+1)
}{2}$ whereas the moduli space of curves of genus $g\geq 2$ has dimension
$3g-3$. As soon as $g\geq 4$, we have $\frac{ g(g+1) }{2} > 3g-3$ and
so not all ppav of dimension $g$ are jacobians. The question then becomes
how else can one parametrize ppav? The first step of a generalization
of the notion of jacobian is the construction of a Prym variety which
we describe below.

\section{Pryms}

Suppose that $C$ is a smooth, complete and irreducible curve of genus
$g+1$ with an \'etale double cover $\pi :\widetilde{C}\rightarrow
C$. Then $\widetilde{C}$ has genus $\tilde{g} := 2g+1$. Let $\sigma
:\widetilde{C}\rightarrow\widetilde{C}$ be the involution of the cover
$\pi$. The involution $\sigma$ acts on the jacobian $J\widetilde{C}$
and the Prym variety $P$ of $\pi$ is an abelian variety of dimension
$g$ defined by
\[
P := im(\sigma - 1)\subset J\widetilde{C}.
\]
The principal polarization of $J\widetilde{C}$ induces twice a
principal polarization on $P$ which we will denote $\Xi$.

A priori, there are two ways of obtaining subvarieties of a Prym
variety: by projection and intersection. Since $P$ is a quotient of
$J\widetilde{C}$ via $\sigma -1$, we can take the images of
subvarieties of $J\widetilde{C}$ in $P$. For $e=1$, we can take the
images in $P$ of Abel curves in $J\widetilde{C}$. We obtain in this
way embeddings of $\widetilde{C}$ in $P$ whose images are called
Prym-embedded curves. The class of a Prym-embedded curve is
$2\frac{[\Xi]^{ g-1 }}{ (g-1 )!}$. Welters classfied all curves of
class twice the minimal class in a ppav \cite{welters871}. The list is
short but it is not limited to Prym-embedded curves. Therefore the
analogue of Matsusaka's theorem is false for curves representing twice
the minimal class. For $e >1$, using Pontrjagin product, it is easy to
see that the projection in $P$ of the image of $\phi_{{\mathcal L}}$
has class $2^e\frac{[\Xi]^{ g-e }}{ (g-e )!}$.

Secondly, since the Prym variety is also a subvariety of
$J\widetilde{C}$, we can intersect the images of the symmetric powers
of $\widetilde{C}$ with it. Since $\Theta_{\widetilde{C}}$ induces
twice $\Xi$ on $P$, for ${\mathcal L}$ of degree $e$, the intersection
of the image of $\phi_{{\mathcal L}}$ with $P$ has class $2^{\tilde{
g}-e }\frac{ [\Xi ]^{\tilde{g}-e }}{ (\tilde{g}-e )!}$ , provided that
the intersection is proper. There is at least one case in which one
can do better: there are a finite number of theta divisors in
$J\widetilde{C}$ whose intersection with $P$ is $2\Xi $ as a divisor
(see for instance \cite{mumford74}). In particular, we have a nice way
of parametrizing the theta divisor of the ppav $(P ,\Xi )$.

It can easily happen that the intersection $\phi_{{\mathcal L} }
(\widetilde{C}^{ (e) })\cap P$ is not proper. In such a case the
cohomolgy class of the resulting subvariety needs to be determined by
other means. Such subvarieties appear in the work of Recillas
\cite{recillas74}, Donagi \cite{donagi81}, Clemens-Griffiths
\cite{clemensgriffiths72} and Beauville \cite{beauville82}. Following
Beauville, we shall call them special subvarieties. A different way of
defining a special subvariety of a Prym variety which allows one to
compute its cohomology class is as follows \cite{recillas74}, \cite{donagi81},
\cite{clemensgriffiths72} and \cite{beauville82}. Let $g^r_e$ be a {\em
complete} linear system of dimension $r$ and degree $e$ on $C$. Let
$L$ be the corresponding invertible sheaf on $C$ and let ${\mathcal
L}$ be an invertible sheaf on $\widetilde{C}$ whose Norm is $Nm
({\mathcal L} ) = L$ (i.e., if ${\mathcal L}\cong{\mathcal
O}_{\widetilde{C}} (D)$, then $L\cong{\mathcal O}_C (\pi_*
D)$). Consider $g^r_e$ as a subvariety of $C^{ (e) }$, isomorphic to
${\mathbb P}^r$. Assuming that $g^r_e$ contains reduced divisors, the
inverse image of $g^r_e$ in $\widetilde{C}^{ (e) }$ is reduced. It
splits as the union of two connected components whose images in
$J\widetilde{C}$ by $\phi_{{\mathcal L}}$ are contained in $P$ and the
translate $P'$ of $P$ such that $P\cup P'$ is the kernel of the Norm
$Nm : J\widetilde{C}\rightarrow JC$. Therefore, after translating one
of these subvarieties, we obtain two subvarieties of $P$. They are
isomorphic if $e$ is odd but not if $e$ even. They both have
cohomology class $2^{ e-2r -1 }\frac{[\Xi]^{ g-r }}{ (g-r )!}$ at
least when $1\leq e\leq 2g+1$ and $e > 2r$ \cite{beauville82}.

To see that these special subvarieties are indeed the nonproper
intersections that we mentioned above, one needs to note that the
fibers of the map $J\widetilde{C}\rightarrow JC$ are translates of
$P\cup P'$. The special subvarieties are intersections of these fibers
with images of maps $\phi_{{\mathcal L}}$. Note that when the $g^r_e$
is nonspecial, i.e., $h^1 (g^r_e) =0$, the special subvariety is in
fact a proper intersection and Beauville's cohomology class is equal
(as it should be) to the cohomology class of the proper intersection
above. As we shall see below, looking at special subvarieties as such
nonproper intersections allows us to define them for arbitrary ppav.

For $g\leq 5$, all ppav are Prym varieties (in the generalized sense
of Beauville \cite{beauville771}). For $g\geq 5$ Prym
varieties of dimension $g$ depend on the same number of moduli as
curves of genus $g+1$, meaning $3g$ moduli. Therefore, for $g\geq 6$,
a general ppav is not a Prym variety. So we need to find a different
way to parametrize a ppav.

\section{Prym-Tjurin varieties}

Again, one would want to use a construction involving curves. Higher
degree coverings $\widetilde{C}\rightarrow C$ do not yield general
ppav because the dimension of the Prym variety (defined as the
connected component of ${\mathcal O}_{\widetilde{C}}$ of the kernel of
the Norm (or pushforward on divisors) $J\widetilde{C}\rightarrow JC$)
is too high and therefore the families of ppav that one would obtain
are too small, their dimensions being the dimensions of the moduli
spaces for the bottom curves $C$. Looking back at a Prym variety, we
note that it was defined as a special type of abelian subvariety of a
jacobian. An abelian subvariety $A$ of a jacobian $(JX ,\Theta_X)$
such that $\Theta_X$ induces $m$-times a principal polarization
$\Theta$ on $A$ is called a Prym-Tjurin variety. Welters has proved
that all ppav are Prym-Tjurin varieties
\cite{welters871}.

To say that $\Theta_X$ induces $m$-times $\Theta$ on $A$ is equivalent
to saying that the class of the image of an Abel curve in $A$ is
$m$-times the minimal class for curves \cite{welters871}. Here we are
taking the image of an Abel curve by the composition
\[
JX\stackrel{\cong }{\longrightarrow}\hat {JX }\longrightarrow\hat
{A }\stackrel{\cong }{\longrightarrow} A
\]
where $\hat{}$ denotes the dual abelian variety, the first and the
last map are induced by the polarizations $\Theta_X$ and $\Theta$
respectively, and the middle map is the transpose of the embedding of
$A$ in $JX$.

Therefore, finding a structure of Prym-Tjurin variety on a ppav is
equivalent to finding a reduced and irreducible curve $\overline{X}$
in $A$ representing $m$-times the minimal class and such that $A$
embeds in the jacobian of the normalization $X$ of
$\overline{X}$. Given such a structure, we can find subvarieties of
$A$ as in the case of Prym varieties: by projection and
intersection. Since the class of the image of an Abel curve is
$m$-times the minimal class, Pontrjagin product shows that the
projection in $P$ of the image of $\phi_{{\mathcal L}}$ has class
$m^e\frac{[\Theta]^{ g-e }}{ (g-e )!}$. Since $\Theta_X$ induces
$m$-times $\Theta$ on $A$, for ${\mathcal L}$ of degree $e$, the
intersection of the image of $\phi_{{\mathcal L}}$ with $P$ has class
$m^{g_X-e}\frac{[\Theta]^{ g_X-e }}{ (g_X-e )!}$, provided that the
intersection is proper ($g_X$ is the genus of $X$). However, unlike
Prym varieties, it is not clear whether one can find translates of
$\Theta_X$ whose intersection with $A$ is $m$-times a theta divisor.
Kanev has shown that this is possible under a restrictive hypothesis
which we explain below.

Any abelian subvariety of $JX$ is the image of an endomorphism of $JX$
(which is not unique). The datum of an endomorphism of $JX$ is
equivalent to the datum of a correspondence, i.e., a divisor in
$X\times X$, up to addition and subtraction of the fibers of the two
projections. This is best seen as follows. Start with a divisor
$D\subset X\times X$. To $D$ one can associate an endomorphism of $JX$
in the following way
\[
\begin{array}{crcl}
\psi_D : & JX &\longrightarrow & JX\\
& {\mathcal O}_X (E) &\longmapsto &{\mathcal O}_X ( {p_2}_*((p_1^* E)\cdot D))
\end{array}
\]
where $p_1$ and $p_2$ are the two projections $X\times X\rightarrow X$. If $D$
is linearly equivalent to a sum of fibers of $p_1$ and $p_2$, then
$\psi_D$ is the zero endomorphism. If we exchange the roles of $p_1$
and $p_2$ in the above definition then $\psi_D$ is replaced by its
image under the Rosati involution. The correspondence $D$ is said to be
symmetric if there are (not necessarily effective) divisors $a$ and $b$
on $X$ such that $D -D^t$ is linearly equivalent to $p_1^* (a) + p_2^*
(b)$, where $D^t$ is the transpose of $D$, i.e., the image of $D$ under
the involution exchanging the two factors of $X\times X$. So $D$ is
symmetric if and only if $\psi_D$ is fixed by the Rosati involution.
We shall assume that this is the case. This is not restrictive since
any abelian subvariety of $JX$ is always the image of an endomorphism
which is fixed by the Rosati involution (see e.g. \cite{welters871}).

Kanev \cite{kanev87} has shown that if the endomorphism can be
represented by a symmetric fixed-point-free correspondence $D$ (i.e.,
the support of $D$ does not intersect the diagonal of $X\times X$),
then one can find theta divisors $\Theta_X$ such that $\Theta_X |_A =
m\Theta$ as divisors. Furthermore, fixing an invertible sheaf
${\mathcal L}_0\in Pic^{ g_X -1 } X$, an invertible sheaf ${\mathcal
L}\in P\subset JX$ is on $\Theta$ if and only if $h^0 ({\mathcal
L}\otimes {\mathcal L}_0)\geq m$ and ${\mathcal L}\not\in\Theta$ if
and only if $h^0 ({\mathcal L} \otimes {\mathcal L}_0)= 0$. This gives
a nice parametrization of $\Theta$ and even allows one to analyze
the singularities of $\Theta$. It is not known however, whether every
ppav is a Prym-Tjurin variety for a (symmetric) fixed-point-free
correspondence. In addition, two correspondences could induce the same
endomorphism of $JX$ while one is fixed-point-free and the other is
not. In general it is difficult to determine whether a given
endomorphism can be induced by a fixed-point-free correspondence.

As we noted above, we can generalize the notion of special
subvarieties to Prym-Tjurin varieties by defining them to be
non-proper intersections of $A$ with images of symmetric powers of
$X$. It would be interesting to compute the cohomology classes of
these special subvarieties and see whether the analogue of Beauville's
formula holds, meaning, the cohomology class of a special subvariety
of dimension $r$ is $m^{ e-2r-1 }\frac{ [\Theta ]^{ g-r }}{ (g-r) !}$.

Welters showed that every principally polarized abelian variety is a
Prym-Tjurin variety \cite{welters871}. Birkenhacke and Lange showed
that every principally polarized abelian variety is a Prym-Tjurin
variety for an integer $m\leq 3^g (g-1)!$ (see
\cite{birkenhakelange92} page 374 Corollary 2.4)\footnote{Their proof
uses $3$-theta divisors. Using the fact that a general $2$-theta
divisor is smooth, the exact same proof would give $m\leq 2^g
(g-1)!$. For abelian varieties with a smooth theta divisor, the same
proof would give $m\leq (g-1)!$. One needs the Lefschetz hyperplane
theorem which also works for mildly singular theta divisors, see
e.g. \cite{goreskymcpherson88} Chapter $2$.}.

\section{Deforming curves}

The question is to find the smallest integer $m$ for which
$m\frac{[\Theta]^{ g-1 }}{ (g-1 )!}$ can be represented by an
algebraic curve. This naturally defines a stratification of the moduli
space ${\mathcal A}_g$ of ppav. Using results of Kanev, Debarre shows
that if $(A ,\Theta)$ is the Prym-Tjurin variety for a symmetric
fixed-point-free correspondence, then either $Sing (\Theta)$ is empty
or its dimension is at least $g- 2m-2$. Since the theta divisor of a
general ppav is smooth, this suggests that, for a general ppav
$A$, the smallest integer $m$ for which there is a curve of class
$m\frac{[\Theta]^{ g-1 }}{ (g-1 )!}$ in $A$ which in addition gives
$A$ a structure of Prym-Tjurin variety should be at least $\frac{ g-1
}{2}$. It is unlikely however that this bound is effective. Debarre
has proved in
\cite{debarre94} that the smallest integer $m$ for which
$m\frac{[\Theta]^{ g-1 }}{ (g-1 )!}$ is the class of an algebraic
curve is at least $\sqrt{\frac{ g}{ 8}} -\frac{1}{4 }$ if $(A
,\Theta)$ is general.

The difficulty is to produce curves in ppav in nontrivial ways. One
approach that we have considered is to deform curves in
jacobians of curves out of the jacobian locus. More precisely, let
$C$ be a curve of genus $g$ with a $g^1_d$ (a pencil of degree $d$).
Define
\[
X_e ( g^1_d ) :=\{ D_e : \exists D\in C^{ (d-e) } \hbox{ such that }
D_e + D\in g^1_d\}\subset C^{ (e) }
\]
(for the precise scheme-theoretical definition see \cite{I13} when
$e=2$ and \cite{I14} for $e>2$). If $d\geq e+1$, the restriction of a
given morphism $\phi_{\mathcal L}$ to $X_e ( g^1_d )$ is nonconstant
and so we can map $X_e ( g^1_d )$ to $JC$. The cohomology class of the
image of $X_e ( g^1_d )$ in $JC$ is $m$-times the minimal class with
$m= {d-2\choose e-1}$. Given a one-parameter infinitesimal deformation
of the jacobian of $C$ out of the jacobian locus ${\mathcal J}_g$ we
ask when the curve $X_e ( g^1_d )$ deforms with it. Infinitesimal
deformations of $JC$ are parametrized by $H^1 (T_{ JC })$ where $T_{
JC }$ is the tangent sheaf of $JC$. The principal polarization
$\Theta_C$ provides an isomorphism between $H^1 (T_{ JC })$ and the
second tensor power $H^1 ({\mathcal O}_C)^{\otimes 2}$. Under this
isomorphism the globally unobstructed deformations of the pair $(JC
,\Theta_C)$ are identified with the symmetric square $S^2 H^1
({\mathcal O}_C)$. Therefore any quadric in the canonical space of $C$
defines a linear form on the space of these infinitesimal
deformations. When we say that an infinitesimal deformation $\eta\in
S^2 H^1 ({\mathcal O}_C)$ is in the annihilator of a quadric, we mean
that it is in the kernel of the corresponding linear form. We prove
the following in \cite{I13}
\vskip20pt
{\bf Theorem}
{\it Suppose $C$ nonhyperelliptic and $d\geq 4$. If the curve $X_2 ( g^1_d
)$ deforms out of ${\mathcal J}_g$ then
\begin{enumerate}
\item either $d = 4$
\item or $d=5$, $h^0 (g^1_5) =
3$ and $C$ has genus $5$ or genus $4$ and only one $g^1_3$.
\end{enumerate}
In the case $g=5$ if $X_2 ( g^1_5 )$ deforms in a direction $\eta\in
S^2 H^1 ({\mathcal C}_C)$ out of ${\mathcal J}_5$, then $\eta$ is in
the annihilator of the quadric $\cup_{ D\in g^1_5 }\langle D\rangle$.}
\vskip20pt

Here $\langle D\rangle$ denotes the span of the divisor $D$ in the
canonical space of $C$. For $d=3$ the image of $X_2 (g^1_3 )$ in $JC$
is an Abel curve and so by the result of \cite{matsusaka59}, the curve
$X_2 (g^1_3 )$ cannot deform out of ${\mathcal J}_g$. For $d=4$, it
follows from the theory of Prym varieties that $X_2 ( g^1_4 )$ deforms
out of ${\mathcal J}_g$ (into the locus of Prym varieties): in fact
$X_2 ( g^1_4 )$ is a Prym-embedded curve
\cite{recillas74}. For $d=5$, $h^0 (g^1_d) = 3$ and $g=4$ (with only
one $g^1_3$) or $g=5$ we believe that $X_2 ( g^1_5 )$ deforms out of
${\mathcal J}_g$ but we do not have a proof of this. An interesting
question is what are these deformations of $( JC,\Theta_C )$ into
which $X_2 ( g^1_5 )$ deforms. Can one describe them in a concrete
geometric way.

For $e > 2$, the analogous result would be the following. The curve
$X_e (g^1_d )$ deforms out of ${\mathcal J}_g$ if and only if
\begin{itemize}
\item either $e= h^0 (g^1_d)$ and $d=2e$
\item or $e= h^0 (g^1_d ) -1$ and $d=2e+1$.
\end{itemize}
We expect this to be true most of the time. There could, however, be
special pairs $(C, g^1_d )$ for which the curve $X_e ( g^1_d )$
deforms out of the jacobian locus but $g^1_d$ does not verify the
above conditions. For instance, so far my calculations \cite{I14} seem
to indicate that if there is a divisor $D\in X$ with $h^0 (D )\geq 2$,
then $X$ might deform in directions $\eta$ whose images in the
projectivization ${\mathbb P} (S^2 H^1 ({\mathcal O}_C))$ are in the span of
the image of $\cup_{ D'\in |D| }\langle D'\rangle$. Finally, note that a
standard Brill-Noether calculation shows that for general curves of
genus $\geq 7$, the smallest $d$ for which they can have $g^1_d$'s
satisfying $d = 2 h^0 (g^1_d)$ is $d= 2g-4$. In such a case the class
of $X_e ( g^1_d )$ is $m$-times the minimal class with $m=
{2g-6\choose g-3}$ which is then what we would find for a general
ppav. We address the case $e>2$ in \cite{I14}.

\section{The genus}

The cohomology class is one discrete invariant that one can associate
to a curve in a ppav. Another discrete invariant is the genus of the
curve. We refer the reader to the nice paper by Bardelli, Ciliberto
and Verra \cite{bardellicilibertoverra95} for a discussion of this.

%\bibliographystyle{amsplain}
%\bibliography{biblio}
%\end{document}

\providecommand{\bysame}{\leavevmode\hbox to3em{\hrulefill}\thinspace}

\end{document}